\documentclass[11pt]{article}
\newcommand{\D}{\displaystyle}

\begin{document}
\thispagestyle{empty}

\begin{center}
\Large
Two Irrational Numbers That Give the Last Non-Zero Digits of $n!$ and $n^n$.
\end{center}

\vspace*{0.2 in}
 
\begin{flushright}
Gregory P. Dresden \\
Washington \& Lee University \\
Lexington, VA 24450

\end{flushright}

\noindent{\sf {\bf Author's Note:} 
This is a slightly revised version of the article that appeared in print in 
{\em Mathematics Magazine} in October of 2001. 
The original proof of Theorem 2 was incorrect; 
I've fixed that mistake here. My thanks to Antonio M. Oller-Marc\'en and Jos\'e Mara Grau for pointing out to me the error. 

Also, Stan Wagon pointed out in a letter to 
{\em Mathematics Magazine} (February 2002) that the question of the periodicity 
of the last non-zero digit of $n!$ (our Theorem 1) appeared several times in 
{\em Crux Mathematicorum} in the 1990's: see 
v.~18 n.~7 (Sep 1992) page 196 for the statement of the problem,
v.~19 n.~8 (Oct 1993) page 228 for an incorrect solution,
v.~19 n.~9 (Nov 1993) page 260 for Stan Wagon's correct solution,
   and
v.~20 n.~2 (Feb 1994) page 44 for another reference. 

I wrote a sequel to this paper, called
``Three Transcendental Numbers From the Last Non-Zero Digits of $n^n$, $F_n$, and $n!$". It appeared in {\em Mathematics Magazine}, April 2008.}

\vskip .35 in

We begin by looking at the pattern formed from the last (i.e.~unit) digit of 
$n^n$.  Since
$1^1 = 1$, $2^2 = 4$, $3^3 = 27$, $4^4 = 256$, and so on, we can easily calculate the first few numbers in our pattern to be
$1, 4, 7, 6, 5, 6, 3,6\dots $.  We construct a decimal number
$N = 0.d_1d_2d_3\dots d_n\dots$
such that the $n^{\mbox{th}}$ digit $d_n$ of $N$ is the last (i.e.~unit) digit of 
$n^n$; that is, 
$N = 0.14765636\dots$.  
In a recent paper [{\bf 1}], R.~Euler and J.~Sadek showed that this $N$ is a rational number with a period of twenty digits:  
\[
N = 0.\overline{14765636901636567490}.
\]
This is a nice result, and we might well wonder if it can be
extended.  Indeed, Euler and Sadek in [{\bf 1}] recommend looking at the last \underbar{non-zero}
digit of $n!$  (If we just looked at the last digit of $n!$,
we would get a very dull pattern of all $0$\rq s, as $n!$ ends in $0$ for
every $n \geq 5$.)

With this is mind, let\rq s define
$\mbox{lnzd}(A)$
to be the last nonzero digit of the  positive integer $A$; it is easy to see that 
$\mbox{lnzd}(A) = A/10^i \mbox{ mod } 10$, where $10^i$ is the largest power of $10$ that divides $A$.
We wish to investigate not only the pattern formed by $\mbox{lnzd}(n!)$,
but also the pattern formed by $\mbox{lnzd}(n^n)$. 
In accordance with [{\bf 1}],
we define the ``factorial\rq\rq\ number
$F = 0.d_1d_2d_3\dots d_n\dots$
to be the infinite decimal
such that each digit
$d_n = \mbox{lnzd}(n!)$,
and 
we define the ``power\rq\rq\ number
$P = 0.d_1d_2d_3\dots d_n\dots$
to be the infinite decimal
such that each digit
$d_n = \mbox{lnzd}(n^n)$,
and we ask whether these numbers are rational (i.e.~are eventually-repeating
decimals) or irrational.

Although the title of this article gives away the secret, we\rq d like to
point out that
at first glance, 
our ``factorial\rq\rq\  number $F$ exhibits a suprisingly high degree 
of regularity, and a fascinating pattern occurs. 
The first few digits of $F$ are easy to calculate:

\noindent\parbox{1.5in}{\begin{eqnarray*}
1! &=& \underbar{1}\\
2! &=& \underbar{2}\\
3! &=& \underbar{6}\\
4! &=& 2\underbar{4}\\
\end{eqnarray*}}
\parbox{1.5in}{\begin{eqnarray*}
5! &=& 1\underbar{2}0 \\
6! &=& 7\underbar{2}0 \\
7! &=& 50\underbar{4}0 \\
8! &=& 403\underbar{2}0 \\
9! &=& 3628\underbar{8}0 \dots
\end{eqnarray*}}
\parbox{2in}{\begin{eqnarray*}
10! &=& 3628\underbar{8}00 \\
11! &=& 39916\underbar{8}00 \\
12! &=& 479001\underbar{6}00 \\
13! &=& 6227020\underbar{8}00 \\
14! &=& 87178291\underbar{2}00
\end{eqnarray*}}
...

Reading the underlined digits, we have:
\[
F = 0.1264\ 22428\ 88682 \dots
\]
Continuing along this path, we have (to forty-nine decimal places):
\[
F = 0.1264\ 22428\ 88682\ 88682\ 44846\ 44846\ 88682\ 22428\ 22428\ 66264 \dots
\]
It is not hard to show that (after the first four digits) $F$ breaks up into five-digit blocks
of the form $x\ x\ 2x\ x\ 4x$, 
where $x \in \{2,4,6,8\}$, and the $2x$ and $4x$ are 
taken mod $10$. 
Furthermore, if we represent these five-digit blocks by symbols
($\dot2$ for $22428$, $\dot4$ for $44846$, $\dot6$ for $66264$, $\dot8$ for
$88682$, and $\dot1$ for the initial four-digit block of $1264$), we have:
\[
F = 0.\dot1\hskip0.375in \dot2\hskip0.375in \dot8\hskip0.375in \dot8\hskip0.375in       
      \dot4\hskip0.375in \dot4\hskip0.375in \dot8\hskip0.375in \dot2\hskip0.375in 
      \dot2\hskip0.375in \dot6\hskip0.375in \dots
\]
Grouping these symbols into blocks of five and then performing more calculations 
(with the aid of  {\tt Maple})
give us $F$ to $249$ decimal places:
\[
F = 0.\dot1\dot2\dot8\dot8\dot4\ 
\dot4\dot8\dot2\dot2\dot6\ 
\dot2\dot4\dot6\dot6\dot8\ 
\dot4\dot8\dot2\dot2\dot6\ 
\dot4\dot8\dot2\dot2\dot6\ 
\dot8\dot6\dot4\dot4\dot2\ 
\dot2\dot4\dot6\dot6\dot8\ 
\dot6\dot2\dot8\dot8\dot4\ 
\dot2\dot4\dot6\dot6\dot8\ 
\dot2\dot4\dot6\dot6\dot8\ 
\dots
\]
The reader will notice additional patterns in these blocks of five symbols (twenty-five digits).  In fact, such patterns exist for any block of size $5^i$.  However, a pattern is different from a period, and doesn\rq t imply that our
decimal $F$ is rational. Consider the classic example of
$0.1\  01\ 001\ 0001\ 00001\ 000001\ \dots$, which has an obvious pattern but
is obviously irrational.  It turns out that our decimal $F$ is also irrational, as
the following theorem indicates:

{\sc Theorem 1.}
{\it 
Let $F = 0.d_1d_2d_3\dots d_n\dots$
be the infinite decimal
such that each digit
$d_n = \mbox{\rm lnzd}(n!)$. Then, $F$ is irrational.
}

As for our ``power\rq\rq\  number $P$, 
it too might seem to be rational at
first glance.  P is only slightly different from Euler and  Sadek's
rational number N, as seen here:  
\begin{eqnarray*}
N &=& 0.14765\  6369\underbar{0}\  16365\  6749\underbar{0}\ 
        14765\  6369\underbar{0}\  16365\  6749\underbar{0} \dots \\
\mbox{and }\ \ 
P &=& 0.14765\  6369\underbar{1}\  16365\  6749\underbar{6}\ 
        14765\  6369\underbar{9}\  16365\  6749\underbar{6} \dots \\
\end{eqnarray*}
(Again, calculations were performed by {\tt Maple}.) 
Despite this striking similarity between $P$ and $N$, it turns out that $P$, like $F$,
 is irrational:

{\sc Theorem 2.}
{\it 
Let $P = 0.d_1d_2d_3\dots d_n\dots$
be the infinite decimal
such that each digit
$d_n = \mbox{\rm lnzd}(n^n)$. Then, $P$ is irrational.
}

Before we begin with the (slightly technical) proofs, let us pause and see
if we can get a feel for
why these two numbers must be irrational. There is no doubt that both $F$ and $P$ are
highly ``regular\rq\rq, in that both exhibit a lot of repetition. The problem is
that there are \underbar{too many}  patterns in the digits, acting on different scales. Taking $P$, for example, we note that there is an obvious pattern (as
shown by Euler and Sadek in [{\bf 1}]) repeating every $20$ digits with 
$1^1$,
$2^2$, $3^3, \dots, 9^9$ and $11^{11}, 12^{12}, \dots, 19^{19}$, but
this is broken by a similar pattern for $10^{10}, 20^{20}, \dots, 90^{90}$ and
$110^{110} \dots 190^{190}$, which repeats every $200$ digits. This, in turn, is
broken by another pattern repeating every $2000$, and so on. A similar behaviour 
is found for $F$, but in blocks of $5$, $25$, $125$, and so on, as mentioned above.
So, in vague terms, there are always ``new patterns\rq\rq\  starting up in the
digits of $P$ and of $F$, and this is what makes them irrational.

Are there some simple observations that we can make about $P$ and $F$
which might help us to prove our theorems?  To start with, we might
notice that every digit of $F$ (except for the first one) is even.  Can we
prove this?  Yes, and without much difficulty:  

{\sc Lemma 1.}
{\it 
For $n \geq 2$, then $\mbox{\rm lnzd}(n!)$ is in $\{2, 4, 6, 8\}$. 
}

\noindent{\it Proof:} The lemma is certainly true for $n=2,3,4$. For $n \geq 5$, we
note that the prime factorization of $n!$ contains more $2$\rq s than $5$\rq s, and
thus even after taking out all the $10$\rq s in $n!$, the quotient will still be even.
To be precise, the number of $5$\rq s in $n!$ (and thus the number of trailing zeros in its
base-$10$ representation) is 
$\D e_5 = \sum_{i=1}^{\infty} \left[n/5^i\right]$, which is strictly less than
the number of $2$\rq s, 
$\D e_2 = \sum_{i=1}^{\infty} \left[n/2^i\right]$ 
(here, $[\cdot]$ represents the greatest integer function). Hence,
$n!/10^{e_5}$ is an even integer not divisible by
$10$, and so $\mbox{lnzd}(n!) =  n!/10^{e_5} \mbox{ mod } 10$, which 
must be in
$\{2,4,6,8\}$.
This completes the proof.

Another helpful observation is to note that the lnzd function appears to be
multiplicative.  For example,
\begin{eqnarray*}
 &\mbox{lnzd}(12) \cdot  \mbox{lnzd}(53) =  2\cdot 3  &=  6, \\
\mbox{and }\ \ &\mbox{lnzd}(12\cdot 53) = \mbox{lnzd}(636) &= 6.
\end{eqnarray*}
However, we note that at times this ``rule\rq\rq\  fails:
\begin{eqnarray*}
&\mbox{lnzd}(15) \cdot  \mbox{lnzd}(22) =  5\cdot 2  &=  10, \\
\mbox{yet }\ \ &\mbox{lnzd}(15\cdot 22) = \mbox{lnzd}(330) &= 3.
\end{eqnarray*}
So, we can only prove a limited form of multiplicativity, but it is
useful none the less:

{\sc Lemma 2.}
{\it 
Suppose $a,b$ are integers such that 
$\mbox{\rm lnzd}(a) \not= 5$,
$\mbox{\rm lnzd}(b) \not= 5$. 
Then, {\rm lnzd} is multiplicative; that is, 
$\mbox{\rm lnzd}(a\cdot b) = 
\mbox{\rm lnzd}(a) \cdot \mbox{\rm lnzd}(b)$ mod $10$.
}

\noindent{\it Proof:} 
Let $x^\prime$ denote the integer $x$ without its trailing zeros; that is,
$x^\prime = x/10^i$, where $10^i$ is the largest power of $10$ dividing $x$.
(Note that $\mbox{lnzd}(x) = x^\prime$ mod $10$.)
By hypothesis, $a^\prime$
and $b^\prime$ are both $\not= 0$ mod $5$, and so 
$(a\cdot b)^\prime \not= 0$ mod $5$ and so 
$(a\cdot b)^\prime = a^\prime \cdot b^\prime$.
Thus,
\begin{eqnarray*}
\mbox{lnzd}(a\cdot b) &=&
\mbox{lnzd}((a\cdot b)^\prime) =
\mbox{lnzd}(a^\prime \cdot b^\prime) =
a^\prime \cdot b^\prime \mbox{ mod } 10 \\
\ \ \ \ \ \ \ \ &=& (a^\prime \mbox{ mod } 10) \cdot (b^\prime \mbox{ mod } 10 )=
\mbox{lnzd}(a^\prime) \cdot \mbox{lnzd}(b^\prime) =
\mbox{lnzd}(a) \cdot \mbox{lnzd}(b).
\end{eqnarray*}
This completes the proof.

 We are now ready to supply the proof of Theorem 1, in which we show
that $F$ is irrational.  The proof is a little technical, but it relies 
first on assuming that $F$ has a repeating decimal expansion, then on
choosing an appropriate multiple of the period $\lambda_0$ and choosing an
appropriate digit $d$, in order to arrive at a contradiction.

\noindent{\it Proof of Theorem 1:}
We argue by contradiction. Suppose $F$ is rational. 
Then $F$ is eventually periodic; let $\lambda_0$ be the
period (i.e. for every $n$
sufficiently 
large, then $d_n = d_{n + \lambda_0}$).
Write $\lambda_0 = 5^i\cdot K$ such that 
$5 \hspace{-0.055in}\not \hspace{-0.02in}|\, K$
 (we acknowledge that $K$ 
could be $1$) and let $\lambda = 2^i\cdot \lambda_0 = 10^i\cdot K$.
Then, $\mbox{lnzd}(\lambda) = \mbox{lnzd}(K)$,
and since $5 \hspace{-0.055in}\not \hspace{-0.02in}|\, K$, then
$10  \hspace{-0.055in}\not \hspace{-0.02in}|\, K$ and so 
$\mbox{lnzd}(K) = K$ mod $10$. Note also that
$\mbox{lnzd}(2 \lambda) = 2K$ mod $10$.
Choose $M$ sufficiently large so that both of the following are true: 
$\mbox{lnzd}(10^M + \lambda ) = \mbox{lnzd}(\lambda)$ (this can easily be
done by demanding that 
$10^M > \lambda$), and for all $n \geq M$, then $d_n = d_{n + \lambda_0}$,
which of course would then equal $d_{n + \lambda}$.
Finally, let 
$d = \mbox{lnzd}( (10^M - 1)!)$.
By Lemma 1, $d \in \{2,4,6,8\}$, and since
$10^M! = (10^M - 1)! \cdot 10^M$, then $d$ also equals
$\mbox{lnzd}(10^M!)$.

Since $\lambda$ is a multiple of the period $\lambda_0$, then if
we let 
$A = 10^M - 1 + \lambda$ and
$B = 10^M -1 + 2\lambda$, then:
\begin{eqnarray*}
d &=&
 \mbox{lnzd}((10^M-1)!) = \mbox{lnzd}(A!) = \mbox{lnzd}(B!)\\
\mbox{and }\ \ \ 
d &=& 
 \mbox{lnzd}(10^M!) = \mbox{lnzd}((A+1)!) = \mbox{lnzd}((B+1)!)
\end{eqnarray*}
Let\rq s now look at the last two terms in the above equation; it is here
we will find our contradiction.
Note that since $\mbox{lnzd}(A!)=d$, then 
$\mbox{lnzd}(A!)\not= 5$. Also,
since $\mbox{lnzd}(A+1) = \mbox{lnzd}(10^M + \lambda)
= \mbox{lnzd}(\lambda) = K$ mod $10$, we know that 
$\mbox{lnzd}(A+1)\not= 5$. Thus, we can apply Lemma 2 to 
$\mbox{lnzd}(A! \cdot (A+1))$
 to get:
\[
d = \mbox{lnzd}((A+1)!) = \mbox{lnzd}(A!) \cdot \mbox{lnzd}(A+1) = d\cdot K 
\mbox{ mod } 10.
\]
Likewise, working with $B$, we find:
\[
d = \mbox{lnzd}((B+1)!) = \mbox{lnzd}(B!) \cdot \mbox{lnzd}(B+1) = d\cdot 2K 
\mbox{ mod } 10.
\]
Combining these two equations, we get:
\[
d = dK \mbox{ mod } 10 \qquad \qquad d = 2dK \mbox{ mod } 10,
\]
and this becomes 
$d(1-K) = 0 = d(1-2K)$ mod $10$. Since $d$ is even, this implies that
$1-K = 0$ mod $5$ and $1-2K = 0$ mod $5$, which is a contradiction. 
Thus, there can be no period $\lambda_0$ and so $F$ is irrational.
This completes the proof.

We now turn our attention to the ``power\rq\rq\  number P derived from the
last non-zero digits of $n^n$.  This part was more difficult, but a major
step was the discovery that the sequence 
$\mbox{lnzd}(100^{100})$,
$\mbox{lnzd}(200^{200})$,
$\mbox{lnzd}(300^{300}) \dots$
was the same as the sequence 
$\mbox{lnzd}(100^{4})$,
$\mbox{lnzd}(200^{4})$,
$\mbox{lnzd}(300^{4}) \dots$.
This relies not only on the fact that $4 | 100$ but also on the fact that 
$a^b = a^{b+4} \mbox{ mod } 10$ for $b > 0$, as used in the
 following lemma:

{\sc Lemma 3.}
{\it 
Suppose $100 \mid x$. Then, 
$\mbox{\rm lnzd}(x^x) = (\mbox{\rm lnzd } x)^4$ mod $10$.
}

\noindent{\it Proof:} 
As in Lemma 2, let $x^\prime$ denote the integer $x$ without its trailing zeros; that is,
$x^\prime = x/10^i$, where $10^i$ is the largest power of $10$ dividing $x$.
Now,
\begin{eqnarray*}
\mbox{lnzd}(x^x) &=& \mbox{lnzd}( (10^ix^\prime)^{ 10^ix^\prime})\\
&=& \mbox{lnzd}( (10^{i\cdot10^ix^\prime})(x^\prime)^{ 10^i\cdot x^\prime})\\
&=& \mbox{lnzd}( (x^\prime)^{ 10^i \cdot x^\prime}).
\end{eqnarray*}
Since 
$10\hspace{-0.055in}\not|\,x^{\prime}$,
 then $10 \hspace{-0.055in}\not|\,(x^\prime)^{ 10^i \cdot x^\prime}$, and so:
\[
\mbox{lnzd}(x^x) 
=  (x^\prime)^{ 10^i \cdot x^\prime} \mbox{ mod } 10.
\] 
Since $100 \mid x$, then $4 \mid 10^i \cdot x^\prime$, and since $(x^\prime)^n = (x^\prime)^{n+4}$ mod $10$
for every positive $n$, we have:
\begin{eqnarray*}
\mbox{lnzd}(x^x) 
&=&  (x^\prime)^{4} \mbox{ mod } 10 \\
&=&  (\mbox{lnzd } x)^{4} \mbox{ mod } 10.
\end{eqnarray*}
This completes the proof.

With Lemma 3 at our disposal, the proof of Theorem 2 is now fairly
easy.

\noindent{\it Proof of Theorem 2:}
Again, we argue by contradiction. Suppose $P$ 
is rational. Let $\lambda_0$ be the
period,
and choose $j$ sufficiently large such that $10^j > 100(\lambda_0 + 1)!$ 
and such that
$\mbox{lnzd}( (10^j + n \lambda_0)^{ 10^j + n \lambda_0}) = 
\mbox{lnzd}( (10^j)^{ 10^j })$ for every positive $n$.
Choosing $n = 100(\lambda_0 + 1)(\lambda_0 -1)!$, we get:
\[
\mbox{lnzd}( (10^j + 100 (\lambda_0 + 1)!)^{ 10^j + 100( \lambda_0 +1)!}) 
= 
\mbox{lnzd}( (10^j)^{ 10^j }).
\]
We reduce the  left side of the above equation by Lemma 3 and the right side is
obviously $1$, so we have:
\[
(\mbox{lnzd} (10^j + 100(\lambda_0 + 1)!))^{4} \mbox{ mod } 10 = 1,
\]
but since $10^j > 100(\lambda_0 + 1)!$ and 
	$\mbox{lnzd}(100(\lambda_0+1)!) = \mbox{lnzd}((\lambda_0+1)!)$, we can 
	rewrite the above equation as:
\[
(\mbox{lnzd}(\lambda_0 + 1)!)^{4} \mbox{ mod } 10 = 1.
\]
Note that by Lemma 1, the only values of  $\mbox{lnzd}((\lambda_0+1)!)$ are 
$2,4,6,$ and $8$, and raising
these to the fourth power mod $10$ gives us:
\[
6 = 1,
\]
which is a contradiction. Thus, $P$ is irrational.
This completes the proof.

We close by asking the obvious (and very difficult)
 question:  Are $F$ and $P$ algebraic or transcendental? 
  I suspect the latter, but it is only a hunch, and I hope some curious
   reader will continue along this interesting line of study.

\end{document}